# On the Jacobian ideal of a nondegenerate power series

Achim Hennings [*]

**Abstract**: Let $f$ be a nondegenerate power series in several variables. We describe a condition for a polynomial $g$ which implies that the product $gf^k$ by a power of $f$ is not contained in the Jacobian ideal of $f$.

**0 Introduction**

Let $P = \mathbb{C}\{x_1, \ldots, x_n\}$ be the ring of convergent power series and $m$ its maximal ideal. Let $f \in m^2$ be a power series, which is nondegenerate with respect to the Newton polyhedron $\Delta := \Gamma_+(f)$ (in the sense of [Kou, 1.19]) and satisfies $f(\ldots, 0, x_i, 0, \ldots) \neq 0$ ($1 \leq i \leq n$). These assumptions imply that the restriction of $f$ to any coordinate plane of $\mathbb{C}^n$ has an isolated singularity at $0$, and $f_i := x_i f_{x_i}$ ($1 \leq i \leq n$) (with $f_{x_i} = \partial f/\partial x_i$) form a system of parameters.

The usual notations in connection with Newton polyhedra are the following: For $g = \sum_m g_m x^m \in \mathbb{C}[\![x_1, \ldots, x_n]\!]$ and $A \subseteq \mathbb{R}^n$ we write $supp(g) := \{m | g_m \neq 0\}$, $\Gamma_+(g) = $ convex hull of $supp(g) + \mathbb{N}^n$, $g_A = \sum_{m \in A} g_m x^m$. Since $\Delta$ coincides with $\mathbb{R}_+^n$ up to a bounded set (in fact $\Delta + \mathbb{R}_+^n \subseteq \Delta$ is enough), the Newton order of $g$ with respect to $\Delta$ is well-defined by $\nu(g) := \sup\{a \in \mathbb{R}_+ | supp(g) \subseteq a\Delta\}$ (where $0\Delta$ has to be taken as $\mathbb{R}_+^n$).

In the next sections, we have to consider families of divisors on a complex manifold. We put $D_I := \bigcap_{i \in I} D_i$, $D(I) := \bigcup_{i \in I} D_i$ if the family $D_i$, $i \in I$, is locally finite.

Our results concern the ideals $j := (f_{x_1}, \ldots, f_{x_n})$ (Jacobian ideal) and $i := (f_1, \ldots, f_n)$ in $P$.

Let $\delta \subseteq \Delta$ be a compact face of dimension $n - 1 - r \in [0, n-1]$, which does not lie in a coordinate hyperplane, and $\sigma := \mathbb{R}_+ \delta$, $A_\sigma := \mathbb{C}[\sigma \cap \mathbb{Z}^n]$, $K_\sigma := \mathbb{C}[\sigma° \cap \mathbb{Z}^n] \subseteq A_\sigma$.[1] $K_\sigma$ is the canonical module of $A_\sigma$ ([Da, 4.6]).

**Theorem (0.1):** Let $h \in P$ and $g := x_1 \ldots x_n h$.

1) $supp(g) \subseteq n\Delta° \Rightarrow h \in i$.
2) $supp(g) \subseteq (n-r)\delta°$ and $0 \neq [g] \in K_\sigma/(f_{1\delta}, \ldots, f_{1\delta})K_\sigma \Rightarrow$
$$Res_0 \begin{bmatrix} f^r h dx_1 \ldots dx_n \\ f_1, \ldots, f_n \end{bmatrix} \neq 0.$$

The ring $A_\sigma$ will be endowed with the ($\mathbb{Q}-$) grading such that the $f_{i\delta}$ are homogeneous of degree 1. We remark that $n - r$ suitable elements of $f_{1\delta}, \ldots, f_{n\delta}$ are maximally linear independent and form a regular sequence for $A_\sigma$ and $K_\sigma$. It follows that the canonical module of $\bar{A}_\sigma := A_\sigma/(f_{1\delta}, \ldots, f_{n\delta})$ is $\bar{K}_\sigma := K_\sigma/(f_{1\delta}, \ldots, f_{n\delta})K_\sigma$. The degree of the socket (highest component) of $\bar{K}_\sigma$ is $n - r$ (Appendix 2, Lemma (A.2.1)).

---

[*] Universität Siegen, Fakultät IV, Hölderlinstrasse 3, D-57068 Siegen
[1] The interior always refers to the affine hull of the polyhedron.

**Corollary (0.2):**

1) Let $a \in A_\sigma$ homogeneous and $0 \neq [a] \in A_\sigma/(f_{1\delta}, \ldots, f_{n\delta})$. Then $f^r a \notin i$, in particular $f^r \tilde{a} \notin j$ if $a = x_1 \ldots x_n \tilde{a}$.
2) The (induced) Newton order of the socket of $P/i$ is $n - v(x_1 \ldots x_n)$.

The purpose of this article is the proof of theorem (0.1).

Proof of corollary (0.2): Ad 1): Let $b \in K_\sigma$ a homogeneous element such that $[ab] \in \overline{K}_\sigma$ generates the socket. Since $x_1 \ldots x_n | b$ (by assumption on $\delta$), i.e. $b = x_1 \ldots x_n \tilde{b}$, theorem (0.1), 2) asserts for $g := ab$, that $f^r a \tilde{b} \notin i$. As the multiplication $x_1 \ldots x_n : P/j \to P/i$ is well defined (and injective) the supplement follows.

Ad 2). Let $\delta$ be the particular face of $\Delta$ for which $(1, \ldots, 1) \in \sigma^\circ$. Then $[x_1 \ldots x_n] \neq 0$ in $\overline{K}_\sigma$, and there is a homogeneous $h \in A_\sigma$ such that the residue class of $g := x_1 \ldots x_n h$ generates the socket of $\overline{K}_\sigma$. By theorem (0.1), 2) $f^r h \notin i$ but by theorem (0.1), 1) $f^r h m \subseteq i$, and $v(f^r h) = v(f^r_\delta h) = v(f^r_\delta g) - v(x_1 \ldots x_n) = n - v(x_1 \ldots x_n)$. Furthermore for $a \in P$ with $v(a) > n - v(x_1 \ldots x_n)$: $v(ax_1 \ldots x_n) \geq v(a) + v(x_1 \ldots x_n) > n$ and $a \in i$ by theorem (0.1), 1).

**Remark (0.3):** Let $\Omega^n := \Omega^n_{\mathbb{C}^n, 0}$ be filtered by $v(g dx_1 \ldots dx_n) = v(gx_1 \ldots x_n)$. Then the map

$$\Omega^n/j\Omega^n \to P/i, \, g dx_1 \ldots dx_n \mapsto gx_1 \ldots x_n,$$

is injective and strict for the quotient filtrations (cf. e.g. [BGMM, B.1.2.3]). It maps the sockets to one another.

The question about the Newton order of the socket of $\Omega^n/j\Omega^n$ was raised in [BGMM].

# 1 Residues

## 1.1 Definitions

We recall the analytic definition and some properties of the local residue, which we need. Instead of integrating over cycles as in [GH], it is more convenient to integrate $C^\infty$ forms obtained by Fubini's theorem, as is of course well known.

Let $U$ be an $n$-dimensional complex manifold[2] and let $V_1, \ldots, V_n \subseteq U$ be hypersurfaces with $V_1 \cap \ldots \cap V_n = \{x\}$. By a partition of unity one finds $C^\infty$-functions $\rho_i$ with $\rho_i = 1$ near $V_i$ and such that $supp(\rho_1 \ldots \rho_n) \subseteq U$ is compact. The residue of an $n$-form $\varphi \in \Omega^n(U \setminus (V_1 \cup \ldots \cup V_n))$ at $x$ along $V_1, \ldots, V_n$ is then defined by

$$Res_{x, V_1, \ldots, V_n}(\varphi) := \frac{\varepsilon_n}{(2\pi i)^n} \int_U \varphi \wedge \bar{\partial}\rho_1 \wedge \ldots \wedge \bar{\partial}\rho_n, \, \varepsilon_n := (-1)^{n(n-1)/2}.$$

If we replace $\rho_1$ by a similar function $\tilde{\rho}_1$, assuming $supp(\tilde{\rho}_1) \subseteq supp(\rho_1)$ without loss of generality, we have $\alpha := \varphi \wedge (\rho_1 - \tilde{\rho}_1)\bar{\partial}\rho_2 \wedge \ldots \wedge \bar{\partial}\rho_n \in \Gamma_c(U \setminus (V_1 \cup \ldots \cup V_n), \mathcal{E}^{2n-1})$, $\int_U d\alpha = 0 = (-1)^n \int_U \varphi \wedge \bar{\partial}(\rho_1 - \tilde{\rho}_1) \wedge \bar{\partial}\rho_2 \wedge \ldots \wedge \bar{\partial}\rho_n$, and the right hand side of the definition remains unchanged. This shows that the residue is well defined. By looking at $\varphi \wedge \rho_1 \bar{\partial}\rho_2 \wedge ... \wedge \bar{\partial}\rho_n$ we also see that the residue is zero if $\varphi$ is regular along $V_1$.

---

[2] These considerations can be generalized to reduced complex spaces, observe [GH, p. 33].

If $V_1 \cap ... \cap V_n = \{x_1, ..., x_m\}$ is a finite set and $U_i \ni x_i$ are disjoint open neighborhoods, one can choose new $\tilde{\rho}_1, ..., \tilde{\rho}_n$ in such a way that $supp(\tilde{\rho}_1 ... \tilde{\rho}_n)$ is compact in $\bigcup_i U_i$. As the right hand side remains the same, we obtain in this case the sum

$$\sum_{x \in V_1 \cap ... \cap V_n} Res_{x, V_1, ..., V_n}(\varphi).$$

We also use the familiar notation $Res_x \begin{bmatrix} \omega \\ f_1, ..., f_n \end{bmatrix}$ for the residue of $\frac{\omega}{f_1 ... f_n}$ along $V_i = \{f_i = 0\}$, $i = 1, ..., n$.

**Example (1.1):** $U = \{|x_i| < \delta_i\} \subseteq \mathbb{C}^n$, $f_1 = x_1^a, f_2, ..., f_n \in \mathcal{O}(U)$ with the only zero $0 \in U$, $\omega = g dx_1 ... dx_n$. For a smaller $\delta_1$ we may choose $\rho_1, ..., \rho_n$ as $\rho_1(x_1), \rho_i(x_2, ..., x_n)$ ($2 \leq i \leq n$). Then the one-dimensional residue

$$Res_0 \left(\frac{h dx_1}{f_1}\right) \text{ with } h := \frac{g}{f_2 ... f_n}$$

is a holomorphic function of $(x_2, ..., x_n)$ on $V_1 \setminus (V_2 \cup ... \cup V_n)$, and we have

$$(2\pi i)^n Res_0 \begin{bmatrix} \omega \\ f_1, ..., f_n \end{bmatrix} = \varepsilon_n \int_U \frac{\omega \wedge \bar{\partial}\rho_1 \wedge ... \wedge \bar{\partial}\rho_n}{f_1 ... f_n}$$

$$= \varepsilon_{n-1} \int_U \frac{g}{f_1 ... f_n} (dx_1 \wedge \bar{\partial}\rho_1) \wedge dx_2 \wedge ... \wedge dx_n \wedge \bar{\partial}\rho_2 \wedge ... \wedge \bar{\partial}\rho_n$$

$$= \varepsilon_{n-1} \int_{\substack{|x_i|<\delta_i, \\ (i \geq 2)}} \left( \int_{|x_1|<\delta_1} \frac{h dx_1 \wedge \bar{\partial}\rho_1}{f_1} \right) dx_2 \wedge ... \wedge dx_n \wedge \bar{\partial}\rho_2 \wedge ... \wedge \bar{\partial}\rho_n$$

$$= (2\pi i)^n Res_{0, V_1 \cap (V_2, ..., V_n)} \left( Res_0 \left(\frac{h dx_1}{f_1}\right) dx_2 \wedge ... \wedge dx_n \right)$$

Here and in the following $V_1 \cap (V_2, ..., V_n)$ abbreviates $(V_1 \cap V_2, ..., V_1 \cap V_n)$.

In particular, one obtains $(2\pi i)^n$ times the coefficient of $x_1^{a_1 - 1} ... x_n^{a_n - 1}$ in $g$ if $f_i = x_i^{a_i}$.

One can directly show that the given definition of the residue satisfies the transformation rule[3] and therefore coincides with other definitions. The sum of the residues depends holomorphically on $\varphi$ since, for a small deformation, $\rho_1, ..., \rho_n$ may be kept unchanged. The residue theorem states that for $n$ hypersurfaces $V_1, ..., V_n$ on a compact manifold $U$ of dimension $n$ with $V_1 \cap ... \cap V_n$ finite the residue sum of $\varphi \in \Omega^n(U \setminus (V_1 \cup ... \cup V_n))$ is zero. This follows as we may take $\rho_i = 1$. (Cf. [GH, ch. 5])

---

[3] By expanding with a second system of parameters.

## 1.2 Proper modification

We consider again the $n$-manifold $U$ and hypersurfaces with $V_1 \cap \ldots \cap V_n = \{x\}$ and we look for the behavior of the residue under a proper modification at $x$.

So let $M$ be a complex manifold and $\pi: M \to U$ a proper holomorphic map such that $D := \pi^{-1}(x)_{red} \subseteq M$ is a divisor and $\pi | M \setminus D \to U \setminus \{x\}$ is an isomorphism. Let $Z_i := \pi'(V_i)$ be the strict transform of $V_i$ in $M$.

The following formula permits to reduce $n$-dimensional residues to $(n-1)$-dimensional ones if $Z_1 \cap \ldots \hat{Z}_i \ldots \cap Z_n$ contains only smooth points of $D$.

**Theorem (1.2):** Under the assumptions

(1) $D \cap Z_1 \cap \ldots \cap Z_n = \emptyset$ and
(2) $D \cap Z_1 \cap \ldots \hat{Z}_i \ldots \cap Z_n$ is finite for some $i \in [1, n]$

we have for $\varphi \in \Omega^n(U \setminus (V_1 \cup \ldots \cup V_n))$ the formula

$$Res_{x,V_1,\ldots,V_n}(\varphi) = (-1)^{i-1} \sum_{p \in D \cap Z_1 \cap \ldots \hat{Z}_i \ldots \cap Z_n} Res_{p,D,Z_1,\ldots \hat{Z}_i \ldots, Z_n}(\pi^*\varphi).$$

Proof: Let e.g. $i = 1$. From $\pi^*\rho_1 = 1$ near $D + Z_1$ and $\pi^*\rho_j = 1$ near $Z_j$ we conclude:

$$(2\pi i)^n \varepsilon_n Res_{x,V_1,\ldots,V_n}(\varphi) = \int_U \varphi \wedge \bar{\partial} \rho_1 \wedge \ldots \wedge \bar{\partial}\rho_n = \int_M \pi^*\varphi \wedge \bar{\partial}\pi^*\rho_1 \wedge \ldots \wedge \bar{\partial}\pi^*\rho_n$$

$$= (2\pi i)^n \varepsilon_n \sum_{p \in (D+Z_1) \cap Z_{[2,n]}} Res_{p,D+Z_1,Z_2,\ldots,Z_n}(\pi^*\varphi)$$

This is the asserted formula because $Z_1 \cap Z_{[2,n]} = \emptyset$.

## 1.3 A special residue computation

Let $f, g_j: U \to \mathbb{C}$ ($1 \le j \le n$) be holomorphic functions with a zero at $x \in U$ and $V := (f)$, $V_j := (g_j)$ the corresponding hypersurfaces. We assume $V_1 \cap \ldots \cap V_n = \{x\}$.

As in 1.2 let $\pi: M \to U$ be a proper modification of $U$ in $x$ with $D := \pi^{-1}(x)_{red}$ a divisor. Let $Z, Z_j$ be the strict transforms of $V, V_j$, and $\tilde{h} := \pi^*h$ for $h \in \mathcal{O}(U)$.

We set up a list of assumptions, which are tailored for the intended case of application: $f$ as in section 0, $g_j \in \mathbb{C}f_1 + \cdots + \mathbb{C}f_n$, $\pi$ a toric resolution of $f$.

**Assumptions** for theorem (1.3):

(1) $D = \bigcup_{i \in I} D_i$ and $\forall J \subseteq I$ the intersection $D_J$ is empty or smooth and connected of codimension $|J|$.
(2) The valuations along $D_i$ satisfy $v_{D_i}(\tilde{f}) = v_{D_i}(\tilde{g}_j) \, \forall j$, i.e. $C := (\tilde{f}) - Z = (\tilde{g}_j) - Z_j$.
(3) $\forall J \subseteq I, |J| = k \ge 1$, it is assumed:
   a) $D_J \cap Z \subseteq D_J$ is a hypersurface or empty.
   b) $D_J \cap Z_{[k,n]} = \emptyset$.
   c) $D_J \cap Z_{[k+1,n]}$ is finite.
(4) $\forall J \subseteq I, |J| = k \ge 1$, there are $c_k^J, \ldots, c_n^J \in \mathbb{C}$ with $\sum_{j=k}^n c_j^J \frac{\tilde{g}_j}{\tilde{f}} = 1$ on $D_J \setminus Z$.

**Theorem (1.3):** Under these assumptions let $l \in [1, n]$ and $\psi \in \Omega^n(U)$ with $\tilde{\psi} := \pi^*(\psi) \in \Gamma(M, \Omega^n(\log D)(-(n+1-l)C))$, i.e. $\tilde{f}^{l-1}\tilde{\psi}/(\tilde{g}_1 \ldots \tilde{g}_n) \in \Gamma(M, \Omega^n(\log D)(Z_1 + \cdots + Z_n))$. Then we have the formula

$$Res_{x,V_1,\ldots,V_n}\left(\frac{f^{l-1}\psi}{g_1\ldots g_n}\right) =$$

$$(-1)^{l-1} \sum_{\substack{i_1,\ldots,i_l \in I \\ \text{different}}} c_1^{\{i_1\}} \ldots c_{l-1}^{\{i_1,\ldots,i_{l-1}\}} \sum_{p \in D_{i_1,\ldots,i_l} \cap Z_{[l+1,n]}} Res_{p,D_{i_1,\ldots,i_l} \cap (Z_{l+1},\ldots,Z_n)} R_{D_{i_1}\ldots D_{i_l}}\left(\frac{\tilde{\psi}}{\tilde{g}_l\ldots\tilde{g}_n}\right).$$

Here $D_{i_1,\ldots,i_l} \cap (Z_{l+1},\ldots,Z_n)$ denotes the tuple of intersections as before, and $R_{D_{i_1}\ldots D_{i_l}}$ the iterated Poincaré residue along $D_{i_1},\ldots,D_{i_l}$. In the case $l = n$ we set $Res_{p,\emptyset}(\alpha) = \alpha$ for $\alpha \in \mathbb{C}$.

Proof: For $J = \{i_1, \ldots, i_l\}$ and $D_J \cap Z_{[l+1,n]} \neq \emptyset$ we need to have $D_J \cap Z_j \neq D_J$ for $j = l$ and $j \in [l+1, n]$ in order to fulfill (3), b) and (3), c). Furthermore, if we put $D' = \sum_{j \notin \{i_1,\ldots,i_l\}} D_j$,

$$R_{D_{i_1}\ldots D_{i_l}}\left(\frac{\tilde{\psi}}{\tilde{g}_l\ldots\tilde{g}_n}\right) \in \Gamma(D_J, \Omega^{n-l}(\log D' \cap D_J)(Z_l + \cdots + Z_n)).$$

The residue of this form at $p$ is well defined by condition (3), b) c). Therefore, the whole formula is well defined. We prove it by induction on $l$.

$\underline{l = 1}$: We transform the right hand side as in example (1.1)

$$\sum_{i \in I} \sum_{p \in D_i \cap Z_{[2,n]}} Res_{p,D_i \cap (Z_2,\ldots,Z_n)} R_{D_i}\left(\frac{\tilde{\psi}}{\tilde{g}_1\ldots\tilde{g}_n}\right) = \sum_{i \in I} \sum_{p \in D_i \cap Z_{[2,n]}} Res_{p,D_i,Z_2,\ldots,Z_n}\left(\frac{\tilde{\psi}}{\tilde{g}_1\ldots\tilde{g}_n}\right).$$

Since $D_{ij} \cap Z_{[2,n]} = \emptyset$ for $i \neq j$, this is also the right hand side in theorem (1.2) for $\varphi = \frac{\psi}{g_1\ldots g_n}$.

$\underline{l \to l+1}$: By induction hypothesis the formula holds for $l$ and $f\psi$ replacing $\psi$.

For $p \in D_{i_1,\ldots,i_l} \cap Z_{[l+1,n]}$ we obtain by (4) (and the linearity of the Poincaré residue)

$$\sum_{p \in D_{i_1,\ldots,i_l} \cap Z_{[l+1,n]}} Res_{p,D_{i_1,\ldots,i_l} \cap (Z_{l+1},\ldots,Z_n)} R_{D_{i_1}\ldots D_{i_l}}\left(\frac{\tilde{f}\tilde{\psi}}{\tilde{g}_l\ldots\tilde{g}_n}\right) =$$

$$\sum_{p \in D_{i_1,\ldots,i_l} \cap Z_{[l+1,n]}} Res_{p,D_{i_1\ldots i_l} \cap (Z_{l+1},\ldots,Z_n)} R_{D_{i_1}\ldots D_{i_l}}\left(\frac{\sum_{j=l}^n c_j^{\{i_1,\ldots,i_l\}} \tilde{g}_j \tilde{\psi}}{\tilde{g}_l\ldots\tilde{g}_n}\right) =$$

$$\sum_{p \in D_{i_1,\ldots,i_l} \cap Z_{[l+1,n]}} Res_{p,D_{i_1\ldots i_l} \cap (Z_{l+1},\ldots,Z_n)} R_{D_{i_1}\ldots D_{i_l}}\left(\frac{c_l^{\{i_1,\ldots,i_l\}} \tilde{\psi}}{\tilde{g}_{l+1}\ldots\tilde{g}_n}\right),$$

where we have used, that only one term has $n - l$ pole components.

The residue theorem, applied with $J = \{i_1, \ldots, i_l\}$, $\alpha := R_{D_{i_1}\ldots D_{i_l}}\left(\frac{\tilde{\psi}}{\tilde{g}_{l+1}\ldots\tilde{g}_n}\right)$, $D' = \sum_{j \notin J} D_j$, gives in the case $D_J \cap Z_{[l+1,n]} \neq \emptyset$

$$\sum_{p \in D_J} Res_{p,D_J \cap (Z_{l+1}+D',Z_{l+2},\ldots,Z_n)}(\alpha) = 0.$$

By example (1.1) and (3), b) we conclude

$$\sum_{p \in D_J \cap Z_{[l+1,n]}} Res_{p,D_J \cap (Z_{l+1},Z_{l+2},\ldots,Z_n)}(\alpha) =$$
$$-\sum_{i_{l+1} \notin J} \sum_{p \in D_{i_{l+1}} \cap D_J \cap Z_{[l+2,n]}} Res_{p,D_{i_{l+1}} \cap D_J \cap (Z_{l+2},\ldots,Z_n)} R_{D_{i_{l+1}}}(\alpha).$$

This formula also holds if $D_J \cap Z_{[l+1,n]} = \emptyset$. Namely, if $D_{i_{l+1}} \cap D_J \cap Z_{[l+2,n]} \neq \emptyset$ for some $i_{l+1}$, by (3), c) $Z_{l+2}, \ldots, Z_n$ cut out hypersurfaces of $D_{i_{l+1}} \cap D_J$, and by (3), b) $Z_{l+1} \cap D_J$ is a hypersurface in $D_J$ or empty. As a result the residue theorem still applies. The induction step is thereby proved.

**Remark (1.4):** If $\frac{\tilde{f}^{l-1}\tilde{\psi}}{\tilde{g}_1 \cdots \tilde{g}_n}$ has only logarithmic poles along $D(I_0) \subseteq D$ for a subset $I_0 \subseteq I$, it is enough to have (3) and (4) only for $J \subseteq I_0$.

## 2 Resolution of $f$

Let $f$ be a nondegenerate power series as in section 0, which is holomorphic on an open neighborhood $U$ of zero. In order to apply theorem (1.3) for the proof of theorem (0.1), we need the details of a toric resolution of $f$. (Cf. e.g. [AGV], [V].)

### 2.1 Torus embeddings

We use the standard notation about torus embeddings from [TE] and [Da], preferably from [Da] if they differ. In particular, we write:

$M = \mathbb{Z}^n$, $M_\mathbb{R} = M \otimes \mathbb{R} \cong \mathbb{R}^n$ with standard base $e_1, \ldots, e_n$,

$N = Hom(M, \mathbb{Z})$, $N_\mathbb{R} = N \otimes \mathbb{R}$ with dual base $e_1^*, \ldots, e_n^*$,

$\langle l, x \rangle := l(x)$ for $l \in N_\mathbb{R}, x \in M_\mathbb{R}$.

Let $\Sigma_0$ be the fan consisting of the cone $\mathbb{R}_+^n := \mathbb{R}_+ e_1^* + \cdots + \mathbb{R}_+ e_n^* \subseteq N_\mathbb{R}$ and all its boundary cones. Any fan $\Sigma$ with supporting set $|\Sigma| = \mathbb{R}_+^n$ is a subdivision of $\Sigma_0$.

To each cone $\sigma \in \Sigma$ is associated its dual cone $\check{\sigma} := \{x \in M_\mathbb{R} | \langle a, x \rangle \geq 0 \ \forall a \in \sigma\} \subseteq M_\mathbb{R}$ and the affine algebraic variety $Spec \ \mathbb{C}[\check{\sigma} \cap M]$.

<u>Remark 1:</u> Let $\sigma \subseteq N_\mathbb{R}$ be a convex rational polyhedral cone and $\check{\sigma} \subseteq M_\mathbb{R}$ its dual. Then the (inclusion reversing) map

$$\text{faces of } \sigma \to \text{faces of } \check{\sigma}, \tau \mapsto \tau^\perp \cap \check{\sigma},$$

is bijective with inverse $\tau \mapsto \tau^\perp \cap \sigma$. We have $\tau^\perp \cap \check{\sigma} = b^\perp \cap \check{\sigma}$ for $b \in \tau^\circ$ and $\dim \tau^\perp \cap \check{\sigma} = \dim \tau^\perp = n - \dim \tau$.

Proof: $\tau^\perp \cap \check{\sigma}$ is a face of $\check{\sigma}$ since $\check{\sigma} \subseteq \{x | \langle a, x \rangle \geq 0 \ \forall a \in \tau\} = \check{\tau}$. We show $(\tau^\perp \cap \check{\sigma})^\perp \cap \sigma = \tau$. The inclusion $\supseteq$ is evident. Let $a \in \sigma \setminus \tau$. The face $\tau$ of $\sigma$ is cut out by a linear form $x \in M_\mathbb{R}$ on $N_\mathbb{R}$ with $x|\tau = 0, x|\sigma \setminus \tau > 0$. Then $x \in \tau^\perp \cap \check{\sigma}$ and $\langle a, x \rangle > 0$, hence $a \notin (\tau^\perp \cap \check{\sigma})^\perp \cap \sigma$.

The equality $\tau^\perp \cap \check{\sigma} = b^\perp \cap \check{\sigma}$ follows from $b \pm \varepsilon a \in \tau$ for $a \in \tau$ and $\varepsilon > 0$ small. To determine the dimension of $\tau^\perp \cap \check{\sigma}$, we may include $\tau$ in a complete chain of faces between $cospan(\sigma)$ and $\sigma$. In the extremal cases the answer is easy.

<u>Remark 2:</u> Let $\sigma \subseteq N_\mathbb{R}$ be a convex rational polyhedral cone with a vertex and $\tau$ a face. Then $\check{\tau} = \check{\sigma} + \mathbb{R}m = \check{\sigma} - \mathbb{N}m$, where $m \in (\tau^\perp \cap \check{\sigma})^\circ \cap M$.

Proof: By remark 1 we have $\tau = m^\perp \cap \sigma$, and $(\sigma_1 \cap \sigma_2)^\vee = \check{\sigma}_1 + \check{\sigma}_2$ holds in general, because $(\check{\sigma}_1 + \check{\sigma}_2)^\vee = (\check{\sigma}_1)^\vee \cap (\check{\sigma}_2)^\vee = \sigma_1 \cap \sigma_2$.

By remark 2, $\mathbb{C}[\check{\tau} \cap M] = \mathbb{C}[\check{\sigma} \cap M]_{x^m}$, and $X_{\check{\tau}} \subseteq X_{\check{\sigma}}$ is an open embedding. For faces $\tau_1, \tau_2 \subseteq \sigma$ we have $X_{\check{\tau}_1} \cap X_{\check{\tau}_2} = X_{(\tau_1 \cap \tau_2)^\vee} \subseteq X_{\check{\sigma}}$. Namely, if $m_i \in \check{\sigma} \cap M$ with $\check{\tau}_i = \check{\sigma} - \mathbb{N}m_i$, then $(\tau_1 \cap \tau_2)^\vee = \check{\tau}_1 + \check{\tau}_2 = \check{\sigma} - \mathbb{N}m_1 - \mathbb{N}m_2 = \check{\sigma} - \mathbb{N}(m_1 + m_2)$.

Therefore, the relation "$p \sim q \Leftrightarrow p \in X_{\check{\sigma}_1}$ and $q \in X_{\check{\sigma}_2}$ have a common preimage under $X_{\check{\sigma}_1} \leftarrow X_{(\sigma_1 \cap \sigma_2)^\vee} \rightarrow X_{\check{\sigma}_2}$" on the disjoint union of all $X_{\check{\sigma}}, \sigma \in \Sigma$, is an equivalence relation, which allows to glue the $X_{\check{\sigma}}, \sigma \in \Sigma$. The resulting scheme $X_\Sigma$ is the torus embedding or toric variety associated to the fan $\Sigma$. The $X_{\check{\sigma}}, \sigma \in \Sigma$, are open affine subsets with the property $X_{\check{\sigma}_1} \cap X_{\check{\sigma}_2} = X_{(\sigma_1 \cap \sigma_2)^\vee}$, $\mathbb{C}[(\sigma_1 \cap \sigma_2)^\vee \cap M] = \mathbb{C}[\check{\sigma}_1 \cap M]\mathbb{C}[\check{\sigma}_2 \cap M]$, which shows that $X_\Sigma$ is separated.[4]

We determine $X_{\check{\sigma}} \setminus \bigcup_{\tau \subsetneq \sigma} X_{\check{\tau}}$. We have $X_{\check{\sigma}} \setminus X_{\check{\tau}} = V(x^m)$, $m \in \alpha^\circ \cap M$, where $\alpha$ is the face of $\check{\sigma}$ corresponding to $\tau \subseteq \sigma$. Since $cospan(\check{\sigma}) \subseteq \check{\sigma}$ corresponds to $\sigma \subseteq \sigma$, we obtain $X_{\check{\sigma}} \setminus \bigcup_{\tau \subsetneq \sigma} X_{\check{\tau}} = V(\{x^m | m \notin cospan(\check{\sigma})\})$. This is the affine torus $T_\sigma := X_{cospan(\check{\sigma})} \subseteq X_{\check{\sigma}}$. Thus, we get the decomposition $X_\Sigma = \bigcup_{\sigma \in \Sigma} T_\sigma$ into disjoint locally closed tori $T_\sigma$ of dimension $\dim cospan(\check{\sigma}) = n - \dim \sigma$.

Consider again $\tau \subseteq \sigma$ in $\Sigma$, $T_\tau \subseteq X_{\check{\tau}} \subseteq X_{\check{\sigma}}$ and the face $\alpha$ of $\check{\sigma}$ corresponding to $\tau \subseteq \sigma$. Then $X_{\check{\sigma}} \setminus \bigcup_{\tau \not\subseteq \rho \subseteq \sigma} X_{\check{\rho}} = V(\{x^m | m \in \beta^\circ, \beta \not\subseteq \alpha\}) = X_\alpha \subseteq X_{\check{\sigma}}$ is irreducible of dimension $\dim \alpha = \dim cospan(\check{\tau})$. On the other hand $X_{\check{\sigma}} \setminus \bigcup_{\tau \not\subseteq \rho \subseteq \sigma} X_{\check{\rho}} = \bigcup_{\gamma \subseteq \sigma} T_\gamma \setminus (\bigcup_{\tau \not\subseteq \rho \subseteq \sigma} \bigcup_{\gamma \subseteq \rho} T_\gamma) = \bigcup_{\tau \subseteq \gamma \subseteq \sigma} T_\gamma$. This implies

$$\overline{T}_\tau \cap X_{\check{\sigma}} = \bigcup_{\tau \subseteq \gamma \subseteq \sigma} T_\gamma, \quad F_\tau := \bigcup_{\sigma \in \Sigma}(\overline{T}_\tau \cap X_{\check{\sigma}}) = \bigcup_{\tau \subseteq \gamma} T_\gamma.$$

We remark that $F_\tau$ is again a torus embedding, which is associated to the projection of the star $St(\tau) = \{\sigma | \sigma \supseteq \tau\}$ from $N = (\mathbb{R}\tau \cap N) \oplus N'$ to a complement $N'$ of $(\mathbb{R}\tau \cap N)$.

A further implication is

$$F_{\sigma_1} \cap F_{\sigma_2} = \bigcup_{\sigma_1 \cup \sigma_2 \subseteq \gamma} T_\gamma.$$

This is $F_\sigma$ with $\sigma = \bigcap_{\tau \supseteq \sigma_1 \cup \sigma_2} \tau$ the smallest cone with faces $\sigma_1, \sigma_2$ if there is such, and $\emptyset$ otherwise. For simplicial fans $\sigma = \sigma_1 + \sigma_2$ if it exists.

We consider now two fans $\Sigma_1, \Sigma_2$ with supporting set $\mathbb{R}_+^n$, $\Sigma_2$ being a subdivision of $\Sigma_1$. That means, every $\sigma \in \Sigma_1$ is a union of cones in $\Sigma_2$. The inclusions $\sigma_2 \subseteq \sigma_1$, $\sigma_i \in \Sigma_i$, define morphisms $X_{\check{\sigma}_2} \rightarrow X_{\check{\sigma}_1}$. They are compatible: For $\sigma_2, \sigma_2' \in \Sigma_2$ and $\sigma_1, \sigma_1' \in \Sigma_1$ with $\sigma_2^{(\prime)} \subseteq \sigma_1^{(\prime)}$ we have a commuting diagram:

$$\begin{array}{ccc} X_{\check{\sigma}_2} & \rightarrow & X_{\check{\sigma}_1} \\ \cup | & & \cup | \\ X_{(\sigma_2 \cap \sigma_2')^\vee} & \rightarrow & X_{(\sigma_1 \cap \sigma_1')^\vee} \end{array}$$

---

[4] To show that $(\sigma_1 \cap \sigma_2)^\vee \cap M = (\check{\sigma}_1 \cap M) + (\check{\sigma}_2 \cap M)$ we use the following fact: There is a hyperplane $m^\perp$ through $\sigma_1 \cap \sigma_2$ separating $\sigma_1, \sigma_2$. This is easily seen by replacing $\sigma_1$ with $\sigma_1 - \sigma_2$. (Cf. [TE, p.24], [Da, 5.4].)

Because of $|\Sigma_1| = |\Sigma_2|$, the resulting morphism $\pi: X_{\Sigma_2} \to X_{\Sigma_1}$ is proper ([TE, p. 26], [Da, 5.5.6]).

Consider $\sigma \in \Sigma_2$ and let $\tau \in \Sigma_1$ be the cone with $\sigma° \subseteq \tau°$, i.e. the smallest cone with $\sigma \subseteq \tau$. Then

$$cospan(\check{\tau}) = cospan(\check{\sigma}) \cap \check{\tau},$$

as follows from $cospan(\check{\tau}) = \tau^\perp = \tau^\perp \cap \check{\tau} = b^\perp \cap \check{\tau} = \sigma^\perp \cap \check{\tau}$, $b \in \sigma \cap \tau°$, by remark 1. The diagram

$$\mathbb{C}[cospan(\check{\tau}) \cap M] \subseteq \mathbb{C}[cospan(\check{\sigma}) \cap M]$$
$$\uparrow \qquad \qquad \uparrow$$
$$\mathbb{C}[\check{\tau} \cap M] \subseteq \mathbb{C}[\check{\sigma} \cap M]$$

(the columns correspond to the closed embeddings $T_\tau \subseteq X_{\check{\tau}}$ and $T_\sigma \subseteq X_{\check{\sigma}}$) shows that

$$\pi | T_\sigma \to T_\tau$$

is a toric projection (i.e. a group homomorphism isomorphic to $pr_1: (\mathbb{C}^*)^r \times (\mathbb{C}^*)^s \to (\mathbb{C}^*)^r$, $r = \dim \tau^\perp$, $r + s = \dim \sigma^\perp$). In consequence

$$\pi^{-1}(T_\tau) = \bigcup_{\sigma° \subseteq \tau°} T_\sigma, \; \pi^{-1}(F_\tau) = \bigcup_{\gamma \supseteq \tau, \gamma \in \Sigma_1} \bigcup_{\sigma° \subseteq \gamma°} T_\sigma,$$

$$\pi^{-1}(X_{\check{\tau}}) = \bigcup_{\gamma \subseteq \tau, \gamma \in \Sigma_1} \bigcup_{\sigma° \subseteq \gamma°} T_\sigma = \bigcup_{\sigma \subseteq \tau, \sigma \in \Sigma_2} T_\sigma = \bigcup_{\sigma \subseteq \tau} X_{\check{\sigma}}.$$

We apply this in particular to $\Sigma$ and $\Sigma_0$. Then

$$\pi^{-1}(0) = \bigcup_{\sigma° \subseteq (\mathbb{R}^n_+)°} T_\sigma.$$

In the case $\Sigma_0 \setminus \{\mathbb{R}^n_+\} \subseteq \Sigma$ (i.e. if $e_1^*, \ldots, e_n^*$ are the only edges of $\Sigma$ in a proper coordinate plane), then each such $\sigma$ has an edge $\lambda$ with $\lambda° \subseteq (\mathbb{R}^n_+)°$ (assuming $\Sigma_0 \neq \Sigma$), and

$$\pi^{-1}(0) = \bigcup_{\lambda \in L} F_\lambda, \text{ where } L := \{\lambda \in \Sigma | \dim \lambda = 1, \lambda° \subseteq (\mathbb{R}^n_+)°\}.$$

We have $X_\Sigma \setminus F_\lambda = \bigcup_{\sigma \not\supseteq \lambda} X_{\check{\sigma}} = X_{\Sigma \setminus St(\lambda)}$, $St(\lambda) := \{\sigma \in \Sigma | \sigma \supseteq \lambda\}$, hence $X_\Sigma \setminus \pi^{-1}(0) = \bigcap_{\lambda \in L} X_{\Sigma \setminus St(\lambda)} = X_{\Sigma_0 \setminus \{\mathbb{R}^n_+\}} = \mathbb{C}^n \setminus \{0\}$.

For a resolution of $f$ one considers such fans $\Sigma$ with supporting set $\mathbb{R}^n_+$, for which the supporting function $s_\Delta: |\Sigma| \to \mathbb{R}$, $s_\Delta(a) := \min a(\Delta)$, is linear on all $\sigma \in \Sigma$.

There is a coarsest fan $\Sigma_\Delta$ with this property: For each face $\delta$ of $\Delta$ (including $\Delta$ itself) we take an inner point $p_\delta \in \delta°$. Then the dual cones of $\Sigma_\Delta$ are those generated by $\Delta$ at the $p_\delta$. The cones themselves can be defined by $\sigma(\delta) := \{a \in \mathbb{R}^n_+ | s_\Delta(a) = a(p_\delta)\} \subseteq N_\mathbb{R}$. This is a convex rational polyhedral cone as it can also be described by $\min a(E) \geq a(p_\delta)$, where $E$ is the finite set of vertices of $\Delta$. Obviously $\sigma(\delta) = \{a \in \mathbb{R}^n_+ | s_\Delta(a) = a(x) \; \forall \; x \in \delta\}$ and $\bigcup_{\delta \subseteq \Delta} \sigma(\delta) = \mathbb{R}^n_+$.

The dual of $\sigma(\delta)$ is indeed $\mathbb{R}_+(\Delta - p_\delta)$ because $\sigma(\delta)$ is the dual of the cone $\mathbb{R}_+(\Delta - p_\delta)$ by definition.

If $\delta, \delta'$ are faces of $\Delta$, then $\sigma(\delta) \cap \sigma(\delta') = \{a \in \mathbb{R}_+^n | s_\Delta(a) = a(p_\delta) = a(p_{\delta'})\} = \sigma(\delta) \cap (p_{\delta'} - p_\delta)^\perp$ is a face of $\sigma(\delta)$ because $\langle \sigma(\delta), p_{\delta'} - p_\delta \rangle \geq 0$. Let $\gamma$ be the smallest face of $\Delta$, which contains $\delta$ and $\delta'$. Then $s_\Delta(a) = a(p_\delta) = a(p_{\delta'})$ is equivalent to $s_\Delta(a) = a(p_\gamma)$, and therefore $\sigma(\delta) \cap \sigma(\delta') = \sigma(\gamma)$.

We have shown that $\Sigma_\Delta = \{\sigma(\delta) | \delta \subseteq \Delta\}$ is a fan and has the desired properties. As $\mathbb{R}_+^n \setminus \Delta$ is bounded, the proper boundary cones of $\mathbb{R}_+^n$ belong to $\Sigma_\Delta$, i.e. $\Sigma_0 \setminus \{\mathbb{R}_+^n\} \subseteq \Sigma_\Delta$.

## 2.2 Resolution of $f$ with special properties

From now on let $\Sigma$ be a regular subdivision of $\Sigma_\Delta$ with $\Sigma_0 \setminus \{\mathbb{R}_+^n\} \subseteq \Sigma$ (cf. [TE, p.32], [Da, 8.2.3]). Thereby is defined a smooth variety $X \coloneqq X_\Sigma$ and a proper morphism $\pi: X \to \mathbb{C}^n$, which is an isomorphism above $\mathbb{C}^n \setminus \{0\}$.

The exceptional divisors in $X$ are $F_\lambda$, $\lambda \in L$, where

$$L = \{\sigma \in \Sigma | \dim \sigma = 1, \sigma \neq \mathbb{R}_+ e_i^* \ (1 \leq i \leq n)\}$$

as above. Let $l_\lambda$ be the primitive vector in $\lambda$. Let $\pi_U: X_U \to U$ be the restriction to the inverse image of the open set $U$. For a holomorphic function $h \in \mathcal{O}(U)$ we have

$$(\pi_U^*(h)) = \sum_{\lambda \in L} v_\lambda(h) F_\lambda + \pi_U'((h)),$$

where $v_\lambda(h) = \min l_\lambda(supp(h)) = s_{\Gamma_+(h)}(l_\lambda)$ are the multiplicities and $\pi_U'((h))$ is the strict transform. This follows from $F_\lambda \cap X_{\tilde{\lambda}} = T_\lambda = (x^m)$ for any $m \in M$ with $l_\lambda(m) = 1$. We put $v_\lambda \coloneqq v_\lambda(f)$ and $\tilde{l}_\lambda \coloneqq v_\lambda^{-1} l_\lambda$, in such a way that $\tilde{l}_\lambda^{-1}(1) \cap \Delta$ is a face of $\Delta$.

As the modification considered in sections 1.2, 1.3 we wish to take $\pi_U$, and we have to ensure the requirements there.

**Lemma (2.1):** There is a Zariski-open nonempty set $W \subseteq N_\mathbb{R}^n$, such that for all $(w_1, \ldots, w_n) \in W$ the assumptions of theorem (1.3) for $f$ and $g_j \coloneqq w_j(f)$ $(1 \leq j \leq n)$ are valid. Here $N_\mathbb{C} \coloneqq N \otimes \mathbb{C} \supseteq N_\mathbb{R}$ is considered as the vector space of all derivations of $\mathbb{C}[M]$ of degree zero.

Proof: For any compact face $\delta$ of $\Delta$, $\dim \delta = n - k$, the set $N_\mathbb{R}|\delta = \{l|\delta \mid l \in N_\mathbb{R}\}$ is a vector space of dimension $n - k + 1$ (consisting of the affine functions), and $1 \in N_\mathbb{R}|\delta$. The elements define derivatives of $f_\delta$. Since $f$ is nondegenerate, the derivatives $f_{1\delta}, \ldots, f_{n\delta}$ (where $f_{i\delta} = f_{\delta i}$) have no common zero in $(\mathbb{C}^*)^n$.

Let $W$ be the set of all $(w_1, \ldots, w_n) \in N_\mathbb{R}^n$ satisfying

(i) $\Gamma_+(g_j) = \Delta$;
(ii) $w_1, \ldots, w_n$ is a basis of $N_\mathbb{R}$;
(iii) $\sum_{j=k}^n \mathbb{R} w_j | \delta = N_\mathbb{R} | \delta$ for all $k \in [1, n]$, $\delta \subseteq \Delta$ compact of dimension $n - k$, $\delta^\circ \subseteq (\mathbb{R}_+^n)^\circ$.

We now verify the assumptions of theorem (1.3).

Ad (1): $D = \bigcup_{\lambda \in L} F_\lambda \subseteq X \setminus T$, $T = T_{\{0\}}$, is a divisor with normal crossings as $\Sigma$ is regular. The intersections of the $F_\lambda$ are connected by the formula $F_\sigma \cap F_{\sigma'} = F_{\sigma+\sigma'}$ or empty.

Ad (2): This follows from (i).

Ad (3), a): Since $x_i \nmid f$ we have $Z = \overline{\pi^*(f) \cap (X_U \cap T)}$. Let $X_{\breve\sigma} \subseteq X$ be a chart. As $s_\Delta|\sigma$ is linear, the set $\delta := \{x \in \Delta | s_\Delta(a) = a(x) \; \forall a \in \sigma\}$ is a face of $\Delta$ (the largest face with $\sigma \subseteq \sigma(\delta)$). Let $m \in \delta \cap M$. Then $f = x^m f^\sigma$ with $v_\lambda(f^\sigma) = 0 \; \forall \lambda \subseteq \sigma$, $\dim \lambda = 1$, hence $Z \cap (X_{\breve\sigma} \cap X_U) = (f^\sigma)$. Let $T_\sigma \subseteq X_{\breve\sigma}$ be a stratum with $\pi(T_\sigma) = \{0\}$ (i.e. $\sigma^\circ \subseteq (\mathbb{R}_+^n)^\circ$). Then $\delta$ is compact and for every monomial $x^p$ in $x^{-m}(f - f_\delta)$ and some $\lambda \subseteq \sigma$, $\dim \lambda = 1$, we have $v_\lambda(x^p) > 0$. Thus $x^{-m}(f - f_\delta)|T_\sigma = 0$ and $Z \cap T_\sigma = (x^{-m} f_\delta)$, $x^{-m} f_\delta \in \mathbb{C}[cospan(\breve\sigma) \cap M]$. By non-degeneracy, $(x^{-m} f_\delta) \subseteq (\mathbb{C}^*)^n$ is a smooth hypersurface or empty. Since $\mathbb{C}[cospan(\breve\sigma) \cap M] \subseteq \mathbb{C}[M]$ defines a toric projection $(\mathbb{C}^*)^n \to T_\sigma$, the same holds for $Z \cap T_\sigma$.

Ad (3), b): Consider $k \in [1, n]$ and a stratum $T_\sigma \subseteq D$ of dimension $\leq n - k$. Using $\Gamma_+(g_j) = \Delta$, we obtain $Z_j \cap T_\sigma = (x^{-m} g_{j\delta})$ as before. By (iii) $\sum_{j=k}^n \mathbb{R} g_{j\delta} = \sum_{j=1}^n \mathbb{R} f_{j\delta}$, since $\dim \delta \leq \dim cospan(\breve\sigma) = n - \dim \sigma \leq n - k$. Therefore, $\emptyset = V(x^{-m} f_{1\delta}, \ldots, x^{-m} f_{n\delta}) = Z_{[k,n]} \cap T_\sigma$.

Ad (4): Because of $v_\lambda(f) = v_\lambda(g_j) \; \forall \lambda \in L$, the function $\frac{g_j}{f}$ is holomorphic on $X_U \setminus Z$. Consider again $T_\sigma \subseteq D$, $\dim T_\sigma = n - k$. On $T_\sigma \setminus Z$ we have $\frac{g_j}{f} = \frac{g_{j\delta}}{f_\delta}$, and the claim follows from $f_\delta \in \sum_{j=k}^n \mathbb{R} g_{j\delta}$.

Ad (3), c): Let $T_\sigma \subseteq D$, $\dim T_\sigma = n - k$. By (3), b) $(F_\sigma \setminus T_\sigma) \cap Z_{[k+1,n]} = \emptyset$. Therefore $F_\sigma \cap Z_{[k+1,n]} \subseteq T_\sigma$ is compact and affine, thus finite.

## 3 Proof of theorem (0.1)

Let $f$ be as in section 0, $U \subseteq \mathbb{C}^n$ an open neighborhood of 0 where $f$ is holomorphic and $\pi_U: X_U \to U$ the resolution of section 2.2. Let $w_1, \ldots, w_n \in N_\mathbb{R}$ be a basis, $g_j := w_j(f)$, $V_j := (g_j) \subseteq U$, $Z_j := \pi_U'(V_j)$. To abbreviate we write $dx := dx_1 \ldots dx_n$.

Ad 1): We choose $(w_1, \ldots, w_n)$ as in lemma (2.1). We may assume $h \in \mathcal{O}(U)$, as $i$ contains a power of $m$. The differential form $\pi_U^* \frac{dx}{x_1 \ldots x_n}$ has first order poles along $X \setminus T$ (cf. [Da, § 15]). Because of $v_\lambda(x_1 \ldots x_n h) > v_\lambda(g_1 \ldots g_n) \; \forall \lambda \in L$ and $v_\lambda(x_1 \ldots x_n h) > 0 \; \forall \lambda = \mathbb{R}_+ e_j^* \; (1 \leq j \leq n)$ we have

$$\pi_U^* \frac{h \, dx}{g_1 \ldots g_n} \in \Gamma(X_U, \Omega^n(Z_1 + \cdots + Z_n)),$$

i.e. regularity along $D = \bigcup_{\lambda \in L} D_\lambda$, $D_\lambda := F_\lambda$. By theorem (1.2), $Res_{0, V_1, \ldots, V_n}\left(\frac{h \, dx}{g_1 \ldots g_n}\right) = 0$. Since $\{h | supp(x_1 \ldots x_n h) \subseteq n\Delta^\circ\}$ is an ideal, we can conclude $h \in i = (g_1, \ldots, g_n)$.

Ad 2): For the moment we only demand from $w_j$ that $\Gamma_+(g_j) = \Delta$. The condition $supp(f^r g) \subseteq n\Delta \cap (\mathbb{R}_+^n)^\circ$ implies for $\psi := \frac{g dx}{x_1 \ldots x_n}$

$$\pi_U^* \frac{f^r \psi}{g_1 \ldots g_n} \in \Gamma(X_U, \Omega^n(\log D)(Z_1 + \cdots + Z_n)).$$

Here $D_\lambda$ ($\lambda \in L$) is a pole component

$$\Leftrightarrow v_\lambda(g) = (n-r)v_\lambda(f)$$
$$\Rightarrow \min l_\lambda(\delta^\circ \cap M) = \min l_\lambda(\Delta) = s_\Delta(l_\lambda)$$
$$\Rightarrow l_\lambda|\delta = s_\Delta(l_\lambda) \text{ is constant} \Rightarrow \lambda \subseteq \sigma(\delta).$$

Let $\sigma_1, \ldots, \sigma_t$ be those $\sigma \in \Sigma$ with $\sigma \subseteq \sigma(\delta)$ and $\dim \sigma = \dim \sigma(\delta) = r + 1$. Let $I_s$ be the set of edges of $\sigma_s$ ($I_s \subseteq L$ since $\sigma(\delta)^\circ \subseteq (\mathbb{R}_+^n)^\circ$) and $L_0 := \bigcup_{1 \leq s \leq t} I_s$. Let $V \subseteq N_\mathbb{R}$ be the vector space generated by $\sigma(\delta)$ (or $\sigma_s$). For $J \subseteq L_0$ with $\emptyset \neq J \subseteq I_s$ for some $s \in [1, t]$ we put

$$E_J := \bigcap_{i \in J} \tilde{l}_i^{-1}(1) \subseteq M_\mathbb{R}.$$

$N_\mathbb{R}|E_J = \{l|E_J | l \in N_\mathbb{R}\}$ is the vector space of affine functions on $E_J$, and $V|E_J = \{l|E_J | l \in V\}$ is a subspace of dimension $(r+2) - |J|$. (To see this, take $x_0 \in E_J$. Then the $l_j$, $j \in \bar{J} := I_s \setminus J$, are linear independent on $E_J - x_0$. Therefore, $1|E_J$, $l_j|E_J$, $j \in \bar{J}$, are a basis.)

We choose now $w_1, \ldots, w_n$ in a more restrictive way such that:

(i) $\Gamma_+(g_j) = \Delta \ \forall j$;
(ii) $w_1, \ldots, w_n$ is a basis for $N_\mathbb{R}$;
(iii) $w_i \in V \ \forall i \in [1, r+1]$ (hence a basis by (ii));
(iv) $(w_k, \ldots, w_{r+1}, \ldots, w_n)|E_J$ is a basis for $N_\mathbb{R}|E_J \ \forall J \subseteq L_0$, $|J| = k \geq 1$, with $J \subseteq I_s$ for some $s \in [1, t]$;
(v) $w_i|\delta = 1 \ \forall i \in [1, r+1]$.

To fulfill (iii) and (iv), we choose $w_1, \ldots, w_{r+1} \in V$ such that $w_k, \ldots, w_{r+1}$ remain linear independent in $N_\mathbb{R}|E_J \ \forall J$, $|J| = k$, and then add some more vectors. Since (ii)-(iv) hold for a Zariski-open set in $V^{r+1} \times N_\mathbb{R}^{n-r-1}$ and $\ker w \cap \Delta = \emptyset$ for $w \in \sigma(\delta)$, $w \neq 0$, we can assure (i). The $w_i|\delta$ in (v) are nonzero constant by (iii) and (ii) and can be normalized to get (v).

With this choice, the conditions (1), (2) of theorem (1.3) and (3), (4) for the pole divisor $\bigcup_{i \in L_0} D_i$ instead of $D$ are satisfied. In particular

$$\sum_{j=k}^n c_j^J \frac{\tilde{g}_j}{\tilde{f}} = 1 \text{ on } D_J \setminus Z,$$

if $c_j^J \in \mathbb{R}$ are chosen such that $\sum_{j=k}^n c_j^J w_j|E_J = 1|E_J$. (Cf. proof of lemma (2.1).)

By theorem (1.3) and remark (1.4) we have

$$Res_{0, V_1, \ldots, V_n}\left(\frac{f^r \psi}{g_1 \ldots g_n}\right) =$$

$$(-1)^r \sum_{s=1}^t \sum_{i_1, \ldots, i_{r+1} \in I_s} c_1^{\{i_1\}} \ldots c_r^{\{i_1, \ldots, i_r\}} \sum_{p \in D_{I_s} \cap Z_{[r+2, n]}} Res_{p, D_{I_s} \cap (Z_{r+2}, \ldots, Z_n)} R_{D_{i_1} \ldots D_{i_{r+1}}}\left(\frac{\tilde{\psi}}{\tilde{g}_{r+1} \ldots \tilde{g}_n}\right)$$

.

We investigate one summand in the first sum, e.g. $s = 1$. If we identify $I_1$ with $[1, r+1]$ by a total ordering, this is the residue sum $\sum_{p \in D_{[1,r+1]} \cap Z_{[r+2,n]}} Res_{p, D_{[1,r+1]} \cap (Z_{r+2},...,Z_n)}$ of the differential form

$$c_1 R_{D_1...D_{r+1}} \left( \frac{\tilde{\psi}}{\tilde{g}_{r+1}...\tilde{g}_n} \right).$$

Here

$$c_1 = \sum_{(i_1,...,i_{r+1}) \in P([1,r+1])} c_1^{\{i_1\}} ... c_r^{\{i_1,...,i_r\}} sign(i_1,...,i_{r+1}),$$

with the notation of lemma (3.1) below.

We need a representation of this differential form: Write $\tilde{l}_i = \sum_{j=1}^{r+1} a_{ij} w_j$ ($1 \le i \le r+1$) and $A = (a_{ij}) \in M_{r+1}(\mathbb{R})$. Because of $\tilde{l}_i | \delta = w_i | \delta = 1 \ \forall \ i \in [1, r+1]$ we have $\sum_{j=1}^{r+1} a_{ij} = 1$ ($1 \le i \le r+1$), hence

$$det A = \begin{vmatrix} a_{11} & \cdots & a_{1,r+1} \\ \vdots & \ddots & \vdots \\ a_{r+1,1} & \cdots & a_{r+1,r+1} \end{vmatrix} = \begin{vmatrix} a_{11} & \cdots & a_{1,r} & 1 \\ \vdots & \ddots & \vdots & \vdots \\ a_{r+1,1} & \cdots & a_{r+1,r} & 1 \end{vmatrix},$$

and this is equal to $c_1$ by corollary (3.2) below.

The iterated Poincaré residue $R_{D_1...D_{r+1}} \left( \frac{\tilde{\psi}}{\tilde{g}_{r+1}...\tilde{g}_n} \right)$ can be described as follows: The iterated residue map $R_{D_1...D_{r+1}} = R_{D_{[1,r+1]}} R_{D_{[1,r]}} ... R_{D_1}$ (restricted to $X_{\breve{\sigma}_1}$)

$$\Gamma(X_{\breve{\sigma}_1}, \Omega^n(\log D)) \to \Gamma(T_{\sigma_1}, \Omega^{n-r-1})$$
$$\cap | \qquad \qquad \|$$
$$\mathbb{C}[\breve{\sigma}_1 \cap M] \otimes \Lambda^n M_\mathbb{C} \to \mathbb{C}[cospan(\breve{\sigma}_1) \cap M] \otimes \Lambda^{n-r-1}(V^\perp)_\mathbb{C}$$

is the inner product $i(l_{r+1}) i(l_r) ... i(l_1)$ on the factor $\Lambda^n M_\mathbb{C}$ and the residue map on $\mathbb{C}[\breve{\sigma}_1 \cap M]$.

Therefore we obtain (as $g = g_{(n-r)\delta}$)

$$R_{D_1...D_{r+1}} \left( \frac{\tilde{\psi}}{\tilde{g}_{r+1}...\tilde{g}_n} \right) = \frac{g}{g_{r+1,\delta} ... g_{n,\delta}} i(l_{r+1}) ... i(l_1) \left( \frac{dx}{x} \right)$$

(with $x = x_1 ... x_n$) and

$$i(l_{r+1}) ... i(l_1) = (v_1 ... v_{r+1}) i(\tilde{l}_{r+1}) ... i(\tilde{l}_1) = (v_1 ... v_{r+1}) det(A) i(w_{r+1}) ... i(w_1).$$

In this representation only $(v_1 ... v_{r+1}) det(A)$ depends on $s$. Together with $c_1$ we have the positive factor $(v_1 ... v_{r+1}) det(A)^2$ appearing in the Poincaré residue which depends on $s$.

Because of $cospan(\breve{\sigma}_s) = cospan(\breve{\sigma}(\delta)) \ \forall \ s \in [1, t]$ we may identify all $T_{\sigma_s}$ and obtain the same result for $s \ne 1$ up to a positive factor, as only the positive factor $(v_1 ... v_{r+1}) det(A)^2$ may change.[5]

---

[5] This idea of proof is used in [Vas] in the special case $\dim \delta = 0$. In this case the proof is already finished here.

Because the residue sum is taken over all $p \in D_{I_s} \cap Z_{[r+2,n]} \subseteq T_{\sigma_s}$, also the residue sums for $s = 1, \ldots, t$ differ only by a positive factor.

Therefore we need only to show, that the residue sum is nonzero for $s = 1$. For this purpose, we apply assertion (A.1.1) from the appendix with $X = D_{I_1} = F_{\sigma_1}$ and the divisors $(Z_{r+1}, \ldots, Z_n) \cap D_{I_1}$. We remark that $H^{n-r-1}\Gamma(X, F^{\cdot+1} \otimes \Omega_X^{n-r-1})$ (notation as in (A.1.1)) and the $(n-r)$-th homogeneous component of $K_\sigma / (f_{1\delta}, \ldots, f_{n\delta}) K_\sigma$ (notation as in theorem 0.1) are isomorphic. The residue sum (without the additional factor) is the trace

$$Tr\left[\frac{g}{g_{r+1,\delta} \cdots g_{n,\delta}} i(l_{r+1}) \ldots i(l_1)(\frac{dx}{x})\right]$$

and is nonzero if the class of $g$ in $H^{n-r-1}\Gamma(X, F^{\cdot+1} \otimes \Omega_X^{n-r-1})$ is nonzero. This is the assumption of theorem (0.1), 2), and the theorem is thereby proved.

It remains to supply a lemma and corollary on determinants used in the above proof.

Let $V$ be a real vector space of dimension $n$, $w_1, \ldots, w_n$ a basis for the dual space $V^*$ and $v_1, \ldots, v_{r+1} \in V^*$ linear independent with $v_i = \sum_{j=1}^n a_{ij} w_j$. For $I \subseteq [1, r+1]$, $|I| = k \geq 1$, we consider the set $P(I)$ of all bijective maps $p = (i_1, \ldots, i_k): [1, k] \to I$ and denote by $sign(p)$ the signum of the permutation of $[1, k]$ obtained by identifying $I$ with $[1, k]$ via the natural order. Then $sign(i_1, \ldots, i_k, j) = (-1)^{1+k-s(j)} sign(i_1, \ldots, i_k)$, where $s(j)$ is the position of $j$ in the ordered set $\{i_1, \ldots, i_k, j\}$. For $I \subseteq [1, r+1]$, $J \subseteq [1, n]$, $|I| = |J|$, we write $D(I, J)$ for the corresponding minor of $(a_{ij}) \in M_{r+1,n}(\mathbb{R})$.

**Lemma (3.1):** Assume that for each $I \subseteq [1, r+1]$, $|I| = k \geq 1$, there are coefficients $c_i^I$, $i \in [k, n]$, with $\sum_{i=k}^n c_i^I w_i = 1$ on the affine subspace $E_I := \cap_{j \in I} v_j^{-1}(1)$. Then for all $I \subseteq [1, r+1]$, $|I| = k \geq 1$:

$$\sum_{p=(i_1, \ldots, i_k) \in P(I)} sign(p) c_1^{\{i_1\}} \ldots c_k^{\{i_1, \ldots, i_k\}} = D(I, [1, k]).$$

Proof: First we derive a formula for $c_i^I$. Let $u := \sum_{i=k}^n c_i^I w_i$. From $u|E_I = 1$ we get $u|\cap_{j \in I} \ker v_j = 0$ and $u = \sum_{j \in I} b_j^I v_j$ with $\sum_{j \in I} b_j^I = 1$. By inserting $v_j = \sum_{l=1}^n a_{jl} w_l$ we obtain

$$u = \sum_{i=k}^n c_i^I w_i = \sum_{j \in I} b_j^I \sum_{l=1}^n a_{jl} w_l = \sum_{l=1}^n (\sum_{j \in I} b_j^I a_{jl}) w_l.$$

This amounts to the system of equations

$$\sum_{j \in I} b_j^I = 1,$$

$$\sum_{j \in I} b_j^I a_{jl} = 0, l = 1, \ldots, k-1$$

together with $c_l^I = \sum_{j \in I} b_j^I a_{jl}$, , $l = k, \ldots, n$.

By Cramer's rule

$$b_j^I \left(\sum_{i \in I} (-1)^{1+s(i)} D(I \setminus \{i\}, [1, k-1])\right) = (-1)^{1+s(j)} D(I \setminus \{j\}, [1, k-1]),$$

and inserting this into the formula for $c_l^I$,

$$c_l^I\left(\sum_{i\in I}(-1)^{1+s(i)}D(I\setminus\{i\},[1,k-1])\right) = \sum_{j\in I}(-1)^{1+s(j)}a_{jl}D(I\setminus\{j\},[1,k-1]),$$

where $s(i)$ is the position of $i$ in $I$. For $l = k$ the last formula simplifies to

$$c_k^I\left(\sum_{i\in I}(-1)^{1+s(i)}D(I\setminus\{i\},[1,k-1])\right) = (-1)^{k-1}D(I,[1,k]).$$

We prove now the lemma by induction on $k$.

$k = 1$: $\qquad c_1^{\{i_1\}} = D(\{i_1\},\{1\})$.

$k \to k+1$: Making use of the induction hypothesis, we get

$$\sum_{p=(i_1,\ldots,i_{k+1})\in P(I)} sign(p) c_1^{\{i_1\}} \ldots c_{k+1}^{\{i_1,\ldots,i_{k+1}\}} =$$
$$\sum_{j\in I}\left(\sum_{p'=(i_1,\ldots,i_k)\in P(I\setminus\{j\})} sign(p',j) c_1^{\{i_1\}} \ldots c_k^{\{i_1,\ldots,i_k\}}\right) c_{k+1}^I =$$
$$\sum_{j\in I}(-1)^{1+k-s(j)}\left(\sum_{p'=(i_1,\ldots,i_k)\in P(I\setminus\{j\})} sign(p') c_1^{\{i_1\}} \ldots c_k^{\{i_1,\ldots,i_k\}}\right) c_{k+1}^I =$$
$$\sum_{j\in I}(-1)^{1+k-s(j)}(D(I\setminus\{j\},[1,k])c_{k+1}^I = D(I,[1,k+1]).$$

**Corollary (3.2)**: For $I = [1, r+1]$ we have

$$\sum_{p=(i_1,\ldots,i_{r+1})\in P(I)} sign(p) c_1^{\{i_1\}} \ldots c_r^{\{i_1,\ldots,i_r\}} = \sum_{i\in I}(-1)^{r+1+s(i)}(D(I\setminus\{i\},[1,r]).$$

Proof:

$$\sum_{p=(i_1,\ldots,i_{r+1})\in P(I)} sign(p) c_1^{\{i_1\}} \ldots c_r^{\{i_1,\ldots,i_r\}} =$$
$$\sum_{i\in I}(-1)^{1+r-s(i)}\sum_{p'=(i_1,\ldots,i_r)\in P(I\setminus\{i\})} sign(p') c_1^{\{i_1\}} \ldots c_r^{\{i_1,\ldots,i_r\}} = \sum_{i\in I}(-1)^{1+r-s(i)}(D(I\setminus\{i\},[1,r]).$$

**Appendix 1:** $H^n(X,\Omega_X^n)$ **for toric varieties** (of dimension $n$)

Let $M, N$ be $n$-dimensional mutually dual lattices in $M_\mathbb{R}$, $N_\mathbb{R}$ with pairing $\langle .,. \rangle: N \times M \to \mathbb{Z}$, and $\Delta$ an $n$-dimensional compact convex integral polyhedron in $M_\mathbb{R}$ (i.e. the vertices are in $M$). We consider a fan $\Sigma$ on $N_\mathbb{R}$, such that the support function $s_\Delta(a) = \inf\langle a,\Delta\rangle$ is linear on all $\sigma \in \Sigma$, i.e. $\exists\, m_\sigma \in M \cap \Delta$ with $s_\Delta|\sigma = \langle .,m_\sigma\rangle$. In other words, $\Sigma$ subdivides the fan dual to $\Delta$. For simplicity, we assume that $\Sigma$ is regular.

Associated to these data are the toric variety $X = X_\Sigma$, the torus $T = Spec\,\mathbb{C}[M] \subseteq X$, the largest reduced invariant divisor $D = X \setminus T$, and a $T$-invariant divisor $C$ with order function $ord\,\mathcal{O}(C) = s_\Delta$, i.e. $C|X_{\tilde\sigma} = (x^{-m_\sigma})$. Then

$$\Gamma(X,\mathcal{O}(lC)) = L(l\Delta) := \{g \in \mathbb{C}[M] | supp(g) \subseteq l\Delta\}.$$

Let $g_0, \ldots, g_n \in L(\Delta)$ with common Newton polyhedron $\Delta(g_i) = conv(supp(g_i)) = \Delta$. We have a decomposition $(g_i) = Z_i - C$, where $Z_i$ is a Cartier divisor with no component in $D$. The $g_i$ can be chosen in such a way that

$$Z_0 \cap \ldots \cap Z_n = \emptyset,\ D \cap Z_0 \cap \ldots \hat{Z}_i \ldots \cap Z_n = \emptyset\ \forall\,i \in [0,n].$$

15(This can be deduced from the fact that the linear system $|C|$ has no base points because $(x^{m_\sigma}) + C|X_{\tilde\sigma} = 0$. It also follows from Bertini's theorem as in [Da, 6.8]: For any component $D_i$ of $D$ the line bundle $\mathcal{O}_{D_i}(C)$ on the toric variety $D_i$ is associated to a face of $\Delta$, and the restriction map $\Gamma(X, \mathcal{O}(C)) \to \Gamma(D_i, \mathcal{O}_{D_i}(C))$ is surjective. To see this, we may replace $C$ by an equivalent divisor which does not contain $D_i$. Therefore, for general $g \in L(\Delta)$ the divisor $Z = (g) + C$ intersects all strata $T_\sigma$ of $X$ transversally, and the $g_i$ can be taken as general linear combinations of $g, x_i \partial g / \partial x_i$ ($1 \leq i \leq n$).)

As $T$ contains only finite compact subspaces, the second condition implies

$$Z_0 \cap \ldots \hat{Z}_i \ldots \cap Z_n \text{ is finite.}$$

Let $\mathcal{U}$ be the covering $U_i = X \setminus Z_i$, $i \in [0, n]$, of $X$. It is acyclic for $\Omega_X^n$ as $\Omega_X^n \cong \mathcal{O}(-D)$ and

(1) $H^q(X, \mathcal{O}(-D + lC)) = 0$ for $l > 0, q > 0$ or $l \leq 0, q \neq n$

([Da], [TE]). We also have

(2) $H^0(X, \mathcal{O}(-D + lC)) = L(l\Delta^\circ), l > 0$.

A trace map $Tr: C^n(\mathcal{U}, \Omega_X^n) \to \mathbb{C}$ can be defined by

$$Tr(\varphi) := \sum_{p \in Z_{[1,n]}} Res_{p, Z_1, \ldots, Z_n}(\varphi), \varphi \in \Gamma(X, \Omega_X^n(*(Z_0 + \cdots + Z_n))).$$

By the residue theorem, $Tr(\varphi) = 0$ for $\varphi \in \sum_{i=1}^n \Gamma(X, \Omega_X^n(*(Z_0 + \cdots \hat{Z}_i \ldots + Z_n)))$, and the induced map

$$Tr: H^n(\mathcal{U}, \Omega_X^n) \to \mathbb{C}$$

is well defined.

Consider $\omega := \sum_{i=0}^n (-1)^i \frac{dg_0}{g_0} \wedge \ldots \widehat{\frac{dg_i}{g_i}} \ldots \wedge \frac{dg_n}{g_n}$. We have $\omega \in \Gamma(X, \Omega_X^n(Z_0 + \cdots + Z_n))$. Namely, on the affine chart $X_{\tilde\sigma} \subseteq X$ we have $g_i = x^m \tilde{g}_i$, $Z_i = (\tilde{g}_i)$ with $m = m_\sigma$, $\frac{dg_i}{g_i} = \frac{dx^m}{x^m} + \frac{d\tilde{g}_i}{\tilde{g}_i}$, and

$$\omega = \sum_{i=0}^n (-1)^i \left( \sum_{j<i} (-1)^j \frac{dx^m}{x^m} \frac{d\tilde{g}_0}{\tilde{g}_0} \ldots \widehat{\frac{d\tilde{g}_j}{\tilde{g}_j}} \ldots \widehat{\frac{d\tilde{g}_i}{\tilde{g}_i}} \ldots \frac{d\tilde{g}_n}{\tilde{g}_n} + \right.$$
$$\left. \sum_{j>i} (-1)^{j-1} \frac{dx^m}{x^m} \frac{d\tilde{g}_0}{\tilde{g}_0} \ldots \widehat{\frac{d\tilde{g}_i}{\tilde{g}_i}} \ldots \widehat{\frac{d\tilde{g}_j}{\tilde{g}_j}} \ldots \frac{d\tilde{g}_n}{\tilde{g}_n} + \frac{d\tilde{g}_0}{\tilde{g}_0} \ldots \widehat{\frac{d\tilde{g}_i}{\tilde{g}_i}} \ldots \frac{d\tilde{g}_n}{\tilde{g}_n} \right) = \sum_{i=0}^n (-1)^i \frac{d\tilde{g}_0}{\tilde{g}_0} \ldots \widehat{\frac{d\tilde{g}_i}{\tilde{g}_i}} \ldots \frac{d\tilde{g}_n}{\tilde{g}_n}.$$

For $p \in X_{\tilde\sigma}$ the residue $Res_{p, Z_1, \ldots, Z_n}(\omega) = Res_{p, Z_1, \ldots, Z_n}\left(\frac{d\tilde{g}_1}{\tilde{g}_1} \ldots \frac{d\tilde{g}_n}{\tilde{g}_n}\right)$ is the intersection number $(Z_1 \ldots Z_n)_p$, Therefore $Tr(\omega) = (Z_1 \ldots Z_n)$.

In general, for compact manifolds holds $0 \neq H^n(X, \Omega_X^n) \cong H^{2n}(X, \mathbb{C}) \cong \mathbb{C}$ (by Serre duality). It follows that $Z_1 \cap \ldots \cap Z_n \neq \emptyset$ and

$$Tr: H^n(\mathcal{U}, \Omega_X^n) = H^n(X, \Omega_X^n) \to \mathbb{C}$$

is an isomorphism. In the toric case $(Z_1 \ldots Z_n) = n! Vol(\Delta)$ by [Da, 11.12.2].

1515

A more explicit representation of $H^n(\mathcal{U}, \Omega_X^n)$ is obtained by the Koszul complex defined as follows: $F := \bigoplus_0^n \mathcal{O}(C)$, $F^p := \Lambda^p F$, $s = (g_0, \ldots, g_n) \in \Gamma(X, F)$ and

$$F^\cdot: 0 \to F^0 \to F^1 \xrightarrow{s\wedge} \ldots \to F^{n+1} \to 0.$$

On the chart $X_{\tilde{\sigma}} \subseteq X$ the module $\mathcal{O}(C)$ is free, and $\Gamma(X_{\tilde{\sigma}}, F^\cdot)$ is isomorphic to the Koszul complex of the elements $\tilde{g}_0, \ldots, \tilde{g}_n \in \Gamma(X_{\tilde{\sigma}}, \mathcal{O}_X)$, which is exact because the elements generate the unit ideal.

Because of $F^p = \Lambda^p(\mathcal{O}_X^{n+1} \otimes \mathcal{O}(C)) = \Lambda^p(\mathcal{O}_X^{n+1}) \otimes \mathcal{O}(pC)$ and (1),

$$0 \to \Omega_X^n \to F^1 \otimes \Omega_X^n \to F^2 \otimes \Omega_X^n \to \cdots$$

is an acyclic resolution of $\Omega_X^n$. Therefore

$$H^n(X, \Omega_X^n) = H^n\Gamma(X, F^{\cdot+1} \otimes \Omega_X^n) = \Gamma(X, \Omega_X^n((n+1)C))/\sum_0^n g_i\, \Gamma(X, \Omega_X^n(nC)).$$

The complex $F^\cdot \otimes \Omega_X^n$ is isomorphic to the subcomplex

$$0 \to \Omega_X^n \to \bigoplus_{i=0}^n \mathcal{O}(Z_i) \otimes \Omega_X^n \to \bigoplus_{i<j} \mathcal{O}(Z_i + Z_j) \otimes \Omega_X^n \to \cdots$$

of the Cech resolution

$$0 \to \Omega_X^n \to \bigoplus_{i=0}^n \mathcal{O}(* Z_i) \otimes \Omega_X^n \to \bigoplus_{i<j} \mathcal{O}(* (Z_i + Z_j)) \otimes \Omega_X^n \to \cdots$$

of $\Omega_X^n$. This shows that the map

$$H^n\Gamma(X, F^{\cdot+1} \otimes \Omega_X^n) \to H^n(\mathcal{U}, \Omega_X^n)$$

is an isomorphim. We have proved:

**(A.1.1)** The map $H^n\Gamma(X, F^{\cdot+1} \otimes \Omega_X^n) \to \mathbb{C}$, $[\varphi] \mapsto Tr\left[\dfrac{\varphi}{g_0 \cdots g_n}\right]$ is bijective.

**Appendix 2: Calculation of the socket degree**

Let $\sigma \subseteq \mathbb{R}^n$ be a $k$-dimensional convex rational polyhedral cone with $cospan(\sigma) = \{0\}$. Let $A_\sigma := \mathbb{C}[\sigma \cap \mathbb{Z}^n]$ and $K_\sigma := \mathbb{C}[\sigma^\circ \cap \mathbb{Z}^n]$ be graded by a linear form $\lambda: \mathbb{Q}^n \to \mathbb{Q}$ which is positive on $\sigma \setminus \{0\}$. We consider a homogeneous system of parameters $f_1, \ldots, f_k \in A_\sigma$ with degrees $\alpha_1, \ldots, \alpha_k$.

**(A.2.1) Lemma**: The socket of $\overline{K}_\sigma := K_\sigma/(f_1, \ldots, f_k)K_\sigma$ has degree $\alpha_1 + \cdots + \alpha_k$.

Proof: By taking an integer multiple we can assume $\lambda(\mathbb{Z}^n) \subseteq \mathbb{Z}$. First, let $\sigma = \langle m_1, \ldots, m_k \rangle$, $m_i \in \mathbb{Z}^n$, be simplicial and $f_i = x^{m_i}$. Then $\{[x^m] \mid m \in \mathbb{Z}^n, m = \sum_{i=1}^k a_i\, m_i, 0 < a_i \leq 1\}$ is a $\mathbb{C}$-basis of $\overline{K}_\sigma$, and the socket has degree $\beta := \lambda(m_1) + \cdots + \lambda(m_k)$. The Poincaré series are related by

$$P(K_\sigma) = \frac{P(\overline{K}_\sigma)}{(1 - t^{\lambda(m_1)}) \cdots (1 - t^{\lambda(m_k)})},$$

and $P(K_\sigma)$ is a quotient of two polynomials with degree $\beta$ and leading coefficients 1 resp. $(-1)^k$. The quotient of them gives the value $P(K_\sigma)(\infty)$.

For a general $\sigma$ let $\sigma° = \bigcup_\tau \tau°$ be a decomposition into simplicial open cones (cf. [TE, p. 32], [Da, 8.2]). Then $P(K_\sigma) = \sum_\tau P(K_\tau)$ is also a quotient of polynomials of the same degree and $P(K_\sigma)(\infty) = \sum_\tau P(K_\tau)(\infty) = \sum_\tau (-1)^{\dim \tau}$. This is the negative Euler characteristic $-\chi(B, \partial B) = (-1)^k$ of a $(k-1)$-ball $B$. Again by

$$P(K_\sigma) = \frac{P(\bar{K}_\sigma)}{(1 - t^{\alpha_1}) \ldots (1 - t^{\alpha_k})}$$

we see that $P(\bar{K}_\sigma)$ is a polynomial of degree $\alpha_1 + \cdots + \alpha_k$ with leading coefficient 1.